\documentclass{article}
\usepackage{amsthm,amsfonts, amsbsy, amssymb,amsmath,graphicx}
\usepackage{graphics}

\newtheorem{thm}{Theorem}

\newtheorem{st}{Statement}

\newcounter{tdfn}
 \def\Z{{\mathbb Z}}

\author{Vassily Olegovich Manturov}

\date{}

\title{Free Knots are Not Invertible}

\begin{document}

\maketitle

\section{Introduction. The parity, bracket for free links and the map $\Delta$.}

The goal of the present paper is to show that free knots and links
are in general not invertible: this fact turns out to be
surprisingly non-trivial.

Free links (also known as homotopy classes of Gauss phrases) were
introduced by Turaev \cite{Tur}, and regularly studied by Manturov
\cite{Ma1,Ma2} and Gibson \cite{Gib}.

We first briefly recall the basic definitions from \cite{Ma1}.

By a {\em $4$-graph} we mean a topological space consisting of
finitely many components, each of which is either a circle or a
finite graph with all vertices having valency four.

A $4$-graph is {\em framed} if for each vertex of it, the four
emanating half-edges are split into two sets of edges called {\em
(formally) opposite}.

A {\em unicursal component} of a $4$-graph is either a free loop
component of it or an equivalence class of edges where two edges
$a$,$b$ are called equivalent if there is a sequence of edges
$a=a_{0},\dots, a_{n}=b$ and vertices $v_{1},\dots, v_{n}$ so that
$a_{i}$ and $a_{i+1}$ are opposite at $v_{i+1}$.

As an example of a free graph one may take the graph of a singular
link.

Analogously to $4$-graphs we define {\em long $4$-graphs}; here we
allow two vertices $a,b$ to have valency one (the others having
valency four) in such a way that the edges $x,y$ incident to $a$ and
$b$ should belong to the same unicursal component in the above
sense. One may also think of these two edges of the long $4$-graph
to be {\em noncompact}, i.e., we may think that the vertices of
valency one are removed and the ends of the edges are taken to
infinity.

By a {\em free link} we mean an equivalence class of framed
$4$-valent graphs  modulo the following transformations. For each
transformation we assume that only one fixed fragment  of the graph
is being operated on (this fragment is to be depicted) or some
corresponding fragments of the chord diagram. The remaining part of
the graph or chord diagram are not shown in the picture; the pieces
of the chord diagram not containing chords participating in this
transformation, are depicted by punctured arcs. The parts of the
graph are always shown in a way such that the formal framing
(opposite edge relation) in each vertex coincides with the natural
opposite edge relation taken from  ${\bf R}^{2}$.

\begin{figure}
\centering\includegraphics[width=120pt]{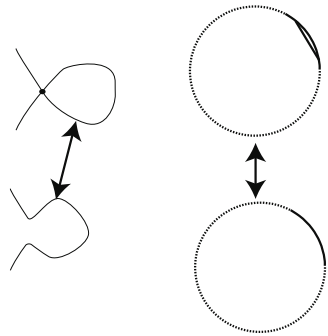}
\caption{Addition/removal of a loop on a graph and on a chord
diagram} \label{1r}
\end{figure}

The first Reidemeister move is an addition/removal of a loop, see
Fig.\ref{1r}

\begin{figure}
\centering\includegraphics[width=150pt]{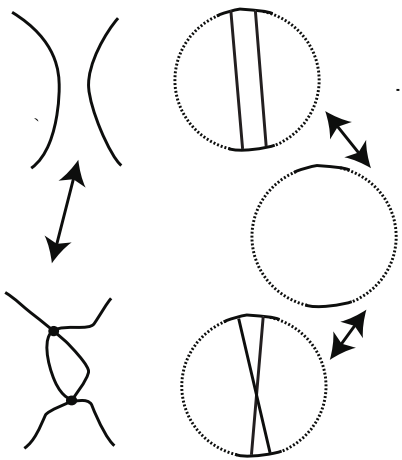} \caption{The second
Reidemeister move and two chord diagram versions of it} \label{2r}
\end{figure}

The second Reidemeister move adds/removes a bigon formed by a pair
of edges which are adjacent in two edges, see Fig. \ref{2r}.

Note that the second Reidemeister move adding two vertices does not
impose any conditions on the edges it is applied to: we may take any
two two edges of the graph an connect them together as shown in Fig.
\ref{2r} to get two new crossings.

\begin{figure}
\centering\includegraphics[width=150pt]{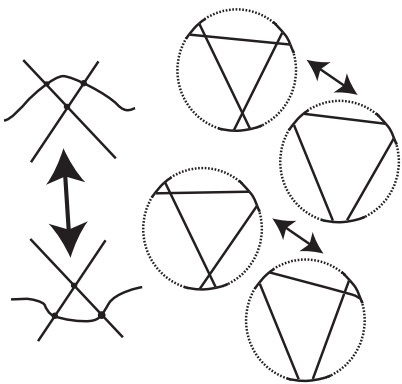} \caption{The third
Reidemeister move and its chord diagram versions} \label{3r}
\end{figure}

The third Reidemeister move is shown in Fig.\ref{3r}.

Note that each of these three moves applied to a framed graph,
preserves the number of unicursal components of the graph. Thus,
applying these moves to graphs with a unique unicursal cycle, we get
to graphs with a unique unicursal cycle.

A {\em free knot} is a free link with one unicursal component
(obviously, the number of unicursal component of a framed $4$-graph
is preserved under Reidemeister moves).

Analogously, one defines {\em long free links and long free knots};
each free link has one {\em noncompact} (long) unicursal component.

Free links are closely connected to  {\em flat virtual knots}, see,
e.g.\cite{EngBook}, i.\,e., with equivalence classes of virtual
knots modulo transformation changing over/undercrossing structure.
The latter are equivalence classes of immersed curves in orientable
$2$-surfaces modulo homotopy and stabilization.

Here we introduce the notion of smoothing, we shall often use in the
sequel.

Let $G$ be a framed four-valent graph, let $v$ be a vertex of $G$
with four incident half-edges $a,b,c,d$, s.t. $a$ is opposite to $c$
and $b$ is opposite to $d$ at $v$.

\begin{figure}
\centering\includegraphics[width=100pt]{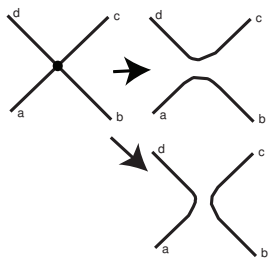} \caption{Two
smoothings of a vertex of for a framed graph} \label{smooth}
\end{figure}

By {\em smoothing} of $G$ at $v$ we mean any of the two framed
$4$-graphs obtained by removing $v$ and repasting the edges as
$(a,b)$, $(c,d)$ or as $(a,d)$ $(b,c)$, see Fig. \ref{smooth}.

Herewith, the rest of the graph (together with all framings at
vertices except $v$) remains unchanged.

We may then consider further smoothings of $G$ at {\em several}
vertices.

We are ready to define invariants of free knots in terms of their
smothings.

\subsection{The parity axioms}

Below we describe an important property of crossings called {\em
parity}, which, if exists, leads to many important consequences in
knot theory: one gets some easy ways for establishing minimality of
knot diagrams, functorial mappings from knots to knots, constructs
powerful invariants, \cite{Ma1,Ma2}.

Assume we have a certain class of knot-like objects which are
equivalence classes of {\bf diagrams} modulo {\em three Reidemeister
moves}. Assume for this class of diagrams (e.g. $4$-valent framed
graphs) there is a fixed rule of distinguishing between two types of
crossings (called even and odd) such that:

1) Each crossing taking part in the first Reidemeister move is even,
and after adding/deleting this crossing the parity of the remaining
crossings remains the same.

2) Each two crossings taking part in the second Reidemeister move
are either both odd or both even, and after performing these moves,
the parity of the remaining crossings remains the same.

3) For the third Reidemeister move, the parities of the crossings
which do not take part in the move remain the same.

\begin{figure}
\centering\includegraphics[width=150pt]{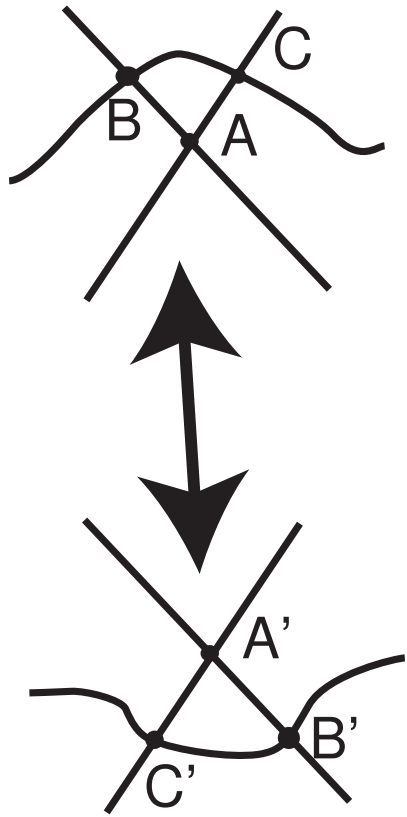} \caption{The third
Reidemeister move} \label{abc}
\end{figure}

Moreover, the parities of the three pairs of crossings are the same
in the following sense: there is a natural one-to-one correspondence
between pairs of crossings $A-A',B-B',C-C'$ taking part in the third
Reidemeister move, see Fig. \ref{abc}.

We require that the {\em parity} of $A$ coincides with that of $A'$,
the {\em parity} of $B$ coincides with that of $B'$ and the parity
of $C$ coincides with that of $C'$.

We also require that the number of odd crossings among the three
crossings in question ($A,B,C$) is even (that is, is equal to $2$ or
$0$).

It turns out that there are {\em many different} parities for
different classes of knots an links.

We should focus on the following two parities:

1) For free knots described by Gauss diagram, we may take a crossing
to be {\em even} if the corresponding chord is linked with an even
number of chords or {\em odd}, otherwise.

2) For free $2$-components links, a crossing formed by one component
is called {\em even}, and a crossing formed by two different
components is called {\em odd}.

The parity axioms in these cases are checked straightforwardly.

\subsection{The bracket $\{\cdot\}$}

\newcommand{\ZGG}{\Z G}

Let $i>2$ be a natural number. Define the set $\ZGG_{i}$ to be the
$\Z_{2}$-linear space generated by set of $i$-component framed
four-valent graphs modulo the following relations:

1) the second Reidemeister moves

2) $L\sqcup \bigcirc=0$, i.\,e., every $n$-component link with one
split trivial component is equivalent to $0$.

For $i=1$, we define $\ZGG_{1}$ analogously with respect to
equivalence 1) and not 2).

It can be easily shown that the elements from any $\ZGG_{i}$ can be
compared algorithmically, namely, each element has a unique minimal
representative which can be obtained by applying consequtively
decreasing second Reidemeister moves.

\newcommand {\K}{\cal K}.

Let $\K$ be the $\Z_2$-linear space generated by free links, and
$\K_{i}$ be the $\Z_{2}$-linear space generated by $i$-component
free links.

 We shall construct a map $\{\cdot\}: \K\mapsto \{K\}$
valued in $\ZGG$ as follows.

Take a framed four-valent graph $G$ representing $K$. By definition,
it has two components. Now, a vertex of $G$ is called {\em odd} if
it is formed by two different components, and {\em even} otherwise.

{\bf The parity axioms can be checked straightforwardly}.

Now, we define

\begin{equation}
\{G\}=\sum_{s} G_{s},
 \label{kbrck}
\end{equation}
where we take the sum over all smoothings of all even vertices, and
consider the smoothed diagrams $K_{s}$ as elements of $\ZGG$. In
particular, we take all elements of $K_{s}$ with free loops to be
zero.

\begin{thm}\cite{Ma1}
The bracket $\{K\}$ is an invariant of two-component free links,
that is, for two graphs $G$ and $G'$ representing the same
two-component free link $K$ we have $\{G\}=\{G'\}$ in $\ZGG$.
\label{mainthm2}
\end{thm}

Analogously one defines $[K]$ as the sum of all one-component
summands (from $\ZGG_{1})$; so, the map $K\to K$ factors through
$\{K\}$.

\subsection{The map $\Delta$}

We shall construct a map from $\K_{2}$ to $\ZGG_{2}$, \cite{Ma2}, as
follows.

In fact, to define the map $\Delta$, one may require for a free knot
to be oriented. However, we can do without.

Given a framed $4$-graph $G$. We shall construct an element
$\Delta(G)$ from $\ZGG_{2}$ as follows. For each crossing $c$ of
$G$, there are two ways of smoothing it. One way gives a knot, and
the other smoothing gives a $2$-component link $G_c$. We take the
one giving a $2$-component link and write

\begin{equation}
\Delta(G)=\sum_{c}G_{c}\in \ZGG_{2}.
\end{equation}

\begin{thm}
$\Delta(G)$ is a well defined mapping from $\K_{1}$ to $\K_{2}$.
\end{thm}

Analogously, one can consider the map $\Delta_{odd}$ where the sum
is taken over all {\em odd crossings} or $\Delta_{even}$ where the
sum is taken over all {\em even crossings}. These maps are both
invariant.

\section{Invertibility of Long Free Knots}

Let us consider the {\em long free knots}. They are defined just as
free knots, but instead of four-valent framed graphs we consider
{\em long four-valent framed graphs}. They can be treated as
four-valent graphs with two infinite edges.

In this section, we consider oriented long free knots. We define the
crossing parity for long free knots just as in the case of the
corresponding compact knots: by using parity of the chords of the
corresponding Gauss diagrams. We are going to prove the following

\begin{thm}

Let $K$ be a framed long four-valent graph with one unicursal
component such that:

1) All crossings of $K$ are odd;

2) There is no room to apply the second decreasing Reidemeister move
to $K$.

3) $K$ is not isomorphic to itself with the orientation reversed.

Then the long free knot represented by $K$ is not invertible.

\label{mindiagg}

\end{thm}

The idea is to modify the bracket for $[K]$ to make it orientable.

Namely, let us define the bracket $[G]_{or}$ for orientable framed
four-valent graphs with one unicursal component as follows. We
define $\ZGG_{1}^{or}$ to be the $\Z_{2}$-linear space of all {\em
oriented} long four-valent framed graphs with one compnonent modulo
the second Reidemeister move.

We take a graph $G$ and take all smoothings of $G$ at even
crossings; each smoothing of such sort is a long framed four-valent
graph; we can naturally endow it with an orientation. Indeed, for
every smoothing $G_{s}$ has two infinite arcs which coincide with
the non-compact arcs of $G$: the initial one $a$ and the final one
$b$. Since $G_{s}$ is a one-component free knot, we may choose these
arcs $a$ and $b$ to be the initial and the final arc, respectively.

Now, analogously to Theorem \ref{mainthm2} one proves the following

\begin{thm}
The bracket $[K]_{or}$ is an invariant of two-component free links,
that is, for two graphs $G$ and $G'$ representing the same
two-component free link $K$ we have $[G]_{or}=[G']_{or}$ in
$\ZGG^{or}_{1}$. \label{mainorient}
\end{thm}

\begin{figure}
\centering\includegraphics[width=200pt]{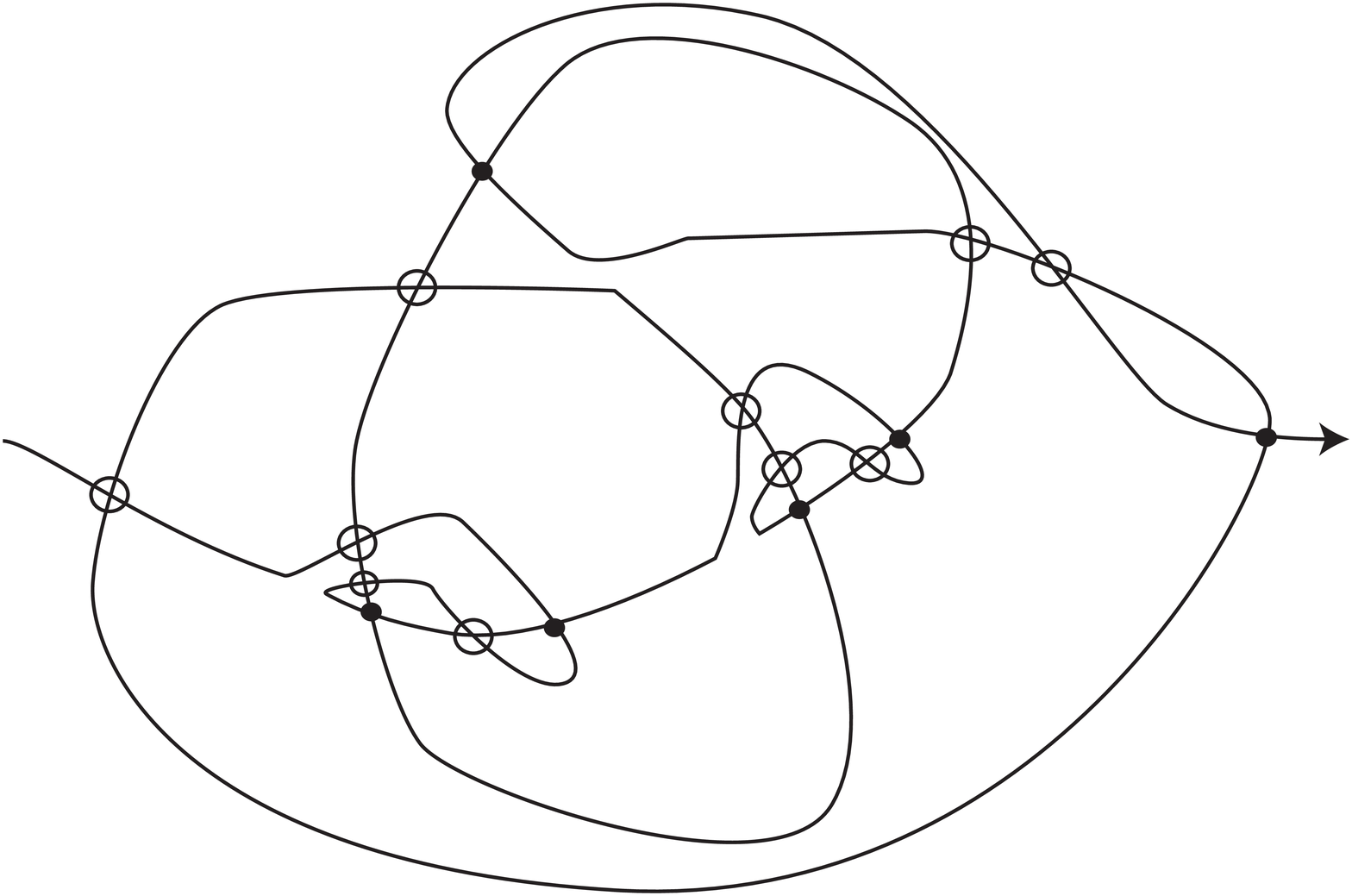} \caption{A
minimal long free knot diagram} \label{minidiagg}
\end{figure}

Now, theorem \ref{mainorient} naturally yields theorem
\ref{minidiagg}. Indeed, if $K$ is a oriented long four-valent free
graph then the equality $[K]_{or}=K$ respects the orientation. Since
$K$ is the minimal representative in its class in $\ZGG^{1}_{or}$,
any other representative $K'$ of the same long free knot has more
crossings. On the other hand, the same is true about ${\bar K}$.
Since $K$ and ${\bar K}$ do not coincide as oriented graphs, the
corresponding knots are different.

Obviously, there are infinitely many examples satisfying theorem
\ref{minidiagg}. One example is shown in Fig. \ref{minidiagg}.

\section{Detecting non-invertibilty of compact links}

As we have seen, the argument of endowing the terms of $[\cdot]$ (or
$\{\cdot\}$) with an orientation works well in the case of long
knots, i.\,e., in the case when we have a reference point. For the
case of compact links (or knots), this is not that easy. Before
defining the ``oriented version'' of the bracket we first collect
the invertibility ``pro'' and ``contra'' arguments in the compact
case.

\subsection{The invertibility arguments}

Let $L$ be an oriented free link, and let ${\bar L}$ be the free
link obtained from $L$ by reversing the orientation of all
components of $L$.

Our goal is to construct such free links $L$ for which ${\bar L}\neq
L$.

Here we collect some observations concerning free knots and links.

\begin{enumerate}

\item The map $\Delta$ can be treated as a map from {\em oriented}
free knots to {\em oriented} two-component free links. Moreover,
${\bar \Delta(K)}=\Delta({\bar K})$. So, having found an example of
a non-invertible free multicomponent link $L$, we may plug in
$\Delta$ in order to get a multicomponent free knot $K$ (trying to
get $\Delta(K)=L+\dots$, where the other summands of $\Delta(K)$ are
immaterial and $L$ yields non-invertibility of $K$.

\item For multicomponent links, one may suggest the following order argument.
For example, consider  a 4-component link $L_{1}\cup L_{2}\cup
L_{3}\cup L_{4}$, where the first component $L_{1}$ has exactly one
intersection point with any other component $L_{i}, i=2,3,4$, whence
any two other components are pairwise disjoint.

We may look at the order of intersection points on the first
component $L_1$ according to its orientation: it can be either
$2,3,4$ or $2,4,3$. Certainly, for the concrete representative
(ordered oriented four-valent framed graph) does not coincide with
its inverse: if $L$ has the order $2,3,4$ then  the inverse link has
the order $2,4,3$.

Nevertheless, these two links are equivalent: the sequence of
Reidemeister moves between these two links is shown in Fig.
\ref{fg}.

\begin{figure}
\centering\includegraphics[width=300pt]{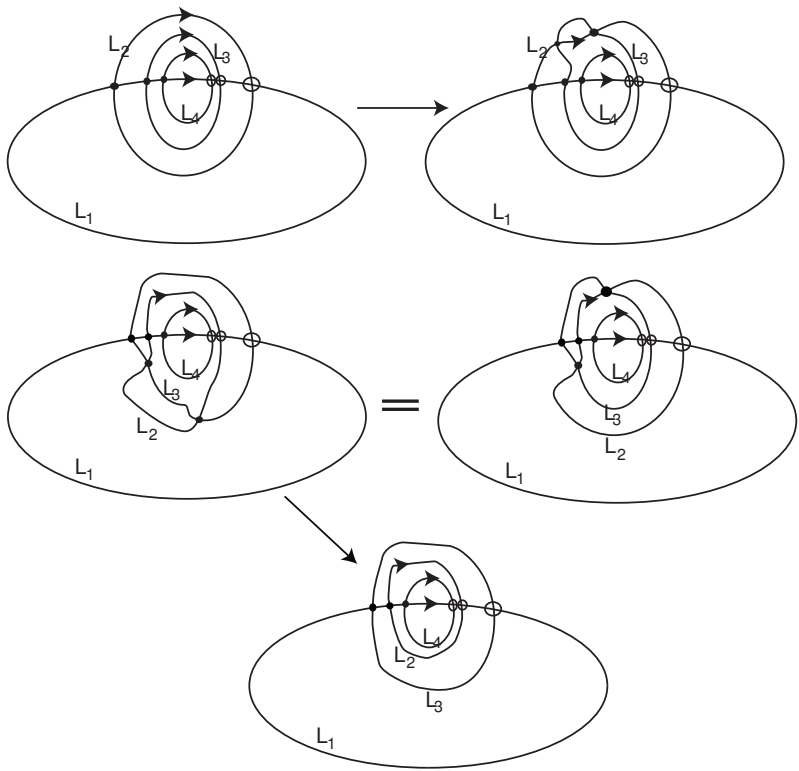} \caption{The
equivalence between a free link and its inverse} \label{fg}
\end{figure}

We first apply the second Reidemeister move to create two
intersection points between the components $L_2$ and $L_3$. Then we
perform a third Reidemeister move for the components $1,2,3$.
Finally, we remove the two intersection points between components
$L_2$ and $L_3$ by a second Reidemeister move.

Finally, we end up with the link where the ordering of the three
intersection points along the orientation of $L_1$ is switched:
instead of $2,3,4$ we get $3,2,4$, which is the just the same as
that for ${\bar L}$.

\item The invariants $[\cdot]$ and $\{\cdot\}$ which sometimes allow one to reduce the
information about a free knot to some information about its
representative graph in the case when $[K]=K$ can not directly be
used for the case of orientable free knots (free links). Indeed, the
bracket $[K]$ (or, in the case of links, its variant $\{L\}$) is
defined as a linear combination of non-oriented four-valent framed
graphs. Indeed, when we try applying third Reiedemeister moves and
collecting terms, we will necessarily get to a situation when a
smoothing at an even crossing breaks the orientation, and we get two
odd crossings where the orientations disagree, \ref{disagree}.

\begin{figure}
\centering\includegraphics[width=200pt]{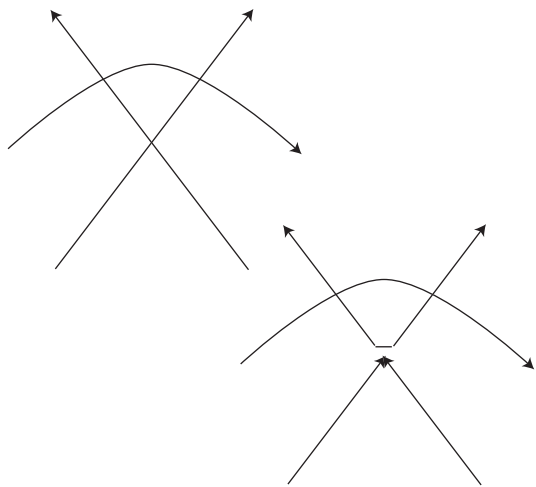}
\caption{Unoriented smoothing} \label{disagree}
\end{figure}

So, the right-hand side in the equality $\{L\}=L$  should be treated
as a non-orientable graph (or, at most, as a partially oriented
graph) modulo second Reidemeister moves. Our goal will be to use
this partial orientability in order to get a genuine orientability.

\end{enumerate}

The main idea of this section is as follows. First, for some
category of two-component links, we modify the bracket $\{L\}$ in
order to make it valued in linear combinations of {\em oriented}
framed graphs modulo second Reidemeister moves compatible with
orientation. This category will include only those two-component
links with orientable atoms. This will lead to some link $L$ where
${\bar L}\neq L$.

Then, by using $\Delta$, we shall extend this result to some
oriented free knots $K$ where $\Delta(K)=L+(\dots)$ $\Delta({\bar
K})={\bar L}+(\dots)$, where the summands $(\dots)$ will mean some
collection of free links which do not affect the non-orientability
of $K$ coming from that of $L$.

Finally, we shall extend this result for knots and links with
orientable atoms.

\subsection{Making the Bracket $\{\cdot\}$ Orientable}

Let ${\cal L}$ be the category of free two-component links
$L_{1}\cup L_{2}$ such that the number of crossing points formed by
both $L_{1}$ and $L_{2}$ is odd.

Obviously, this property is preserved by Reidemeister moves, so, the
category is well-defined.

\newcommand{\Lo}{{\cal L}^{o}}

Now, let $\Lo$ be the category of links from ${\cal L}$ where the
components are ordered: $L_{1}\cup L_{2}$ and the component $L_{1}$
is endowed with an orientation.

The main theorem we are going to prove is the following
\begin{thm}
Let $L$ be a four-valent framed graph having two unicursal
components, one of which is oriented. Assume that:

\begin{enumerate}
\item All crossings of $L$ belong to two different components
$L_{1}$ and $L_{2}$

\item The diagram is irreducible, i.e., no decreasing second
Reidemester move can be applied to it.

\item The diagram is not inveritble, i.e., it is not isomorphic to
itself with the orientation of $L_{1}$ reversed.

\item There is no isomorphism of $L_{1}\cup L_{2}$ onto itself (as framed 4-graphs) which
disregards the orientation and interchanges $L_{1}$ and $L_{2}$.
\end{enumerate}

Then the link $L$ is not invertible.\label{minexm}
\end{thm}

\newcommand{\zLo}{\Z_{2}\Lo}

\newcommand{\tLo}{{\mathfrak{L}}}

Let $\zLo$ be the $\Z_2$-linear space spanned by all links from
$\Lo$. Let $\tLo$ be the quotient linear $\Z_2$-space of the space
spanned by four-valent framed two-component graphs with one
component oriented by the second Reidemeister move.

We want to construct the bracket map $\{\cdot\}:\zLo\to \tLo$ To do
that, we introduce the parity 1.2.1 for links from $\Lo$: a crossing
for a two-component link $L_{1}\cup L_{2}$ is {\em even} if it is
formed by one component $L_{1}$ or $L_{2}$, and it is odd if it is
formed by the two components $L_{1}\cup L_{2}$.

Then for a link $L\in \zLo$ we take its bracket $\{\cdot\}$ as
described above and modify it as follows. First, the bracket
$\{\cdot\}_{2}$ will contain only two-component summands. Note that
every summand of $\{\cdot\}$ has at least two components: some
components correspond to the former $L_{1}$, and the others
correspond to the former $L_{2}$. So, we are interested in the case
when the summand has exactly two components, that is, smoothings at
vertices of $L_{1}$ lead to one component and smoothings at vertices
of $L_{2}$ lead to the the other component.

Note that if we just take $\{\cdot\}_{2}$ to be the sum of these
two-component summands regardless any orientation, it becomes an
invariant of free links, because this map just factors through the
usual $\{\cdot\}$ map.

Now, we would like to endow the summands of $\{\cdot\}_{2}$ with an
orientation of the component $L_{1}$ (by abusing notation we denote
by $L_{1}$ the component consisting of edges belonging to $L_{1}$).

Let $s$ be a smoothing, and let $L_{1}^{s}$ be the result of
applying this smoothing to $L_{1}$ (we agreed that it gives one
component). The number of crossings between the new $L_{1}$ and
$L_{2}$ is the same as that between the old $L_{1}$ and $L_{2}$,
because we do not smooth odd vertices which form crossings between
$L_{1}$ and $L_{2}$.

So, we have some $2n+1$ crossings on the new $L_{1}$ with pieces of
orientation of the original component $L_{1}$ on it. These
orientations may disagree since the way of smoothing of the original
link $L$ does not agree with the orientation of $L_{1}$, in general.

In total, we have $2n+1$ ways of orienting the link $L_{1}$; assume
some $2l+1$ of them give one orienation $o$ and the remaining
$2(n-l)$ ones give the opposite orientation $\bar o$ of $L_{1}$.

Now, we choose the orientation $o$ for $L_{1}$ in the given summand.
So, we endowed the terms of $\{L\}_{2}$ with an orientation of one
component $L_{1}$. From now on we consider $\{L\}_{2}$ as a sum of
two-component framed graphs with one component oriented.

\begin{thm}
The bracket $\{\cdot\}_{2}$ with one-component orientation described
above, is an invariant of two-component links from $\Lo$.
\label{invr}
\end{thm}

\begin{proof}
One should just repeat the invariance proof for the bracket
$\{\cdot\}$ in its usual non-oriented version and see that the
orientation of the $L_{1}$ components for all pairs of cancelling
terms agree.

For the $\Omega_{1}$ move, there is nothing to prove since the only
crossing in question gets smoothed and does not affect the
orientation.

The same happens for the $\Omega_{2}$ move with two crossings on the
same component and for the $\Omega_{3}$ move applied to three
crossings lying on the same component.

Now, for a move $\Omega_{2}$ which is applied to two crossings lying
in $L_{1}$ and $L_{2}$, these two crossings contribute into the
orientation of $L_{1}$; namely, if we had some $2n-1$ crossings
formed by $L_{1}$ and $L_{2}$ before the move, we get $2n+1$
crossings after the move. But the two orientation coming from
initial component $L_{1}$ coming from these two crossings agrees for
the smoothed $L_{1}^{s}$, so the rule for choosing the orientation
for $L_{1}^{s}$ remains the same.

Finally, when we apply the third Reidemeister move referring to two
components, we have to check several cases. If this move applies to
two pieces of $L_{2}$ and one piece of $L_{1}$ then the two
crossings between $L_{2}$ and $L_{1}$ contribute the same
orientation to $L_{1}$ (in the LHS as well as in the RHS, because
they are consequent crossings on the same arc), so the choice of the
orientation remains the same in both sides of the equation.

So, we are left with the case when $L_{1}$ occurs twice and $L_{2}$
occurs once, and the only ``even'' point in our Reidemeister move
belongs to $L_{1}$.

 The two variants of this move are drawn in Fig. \ref{Fg}.

\begin{figure}
\centering\includegraphics[width=200pt]{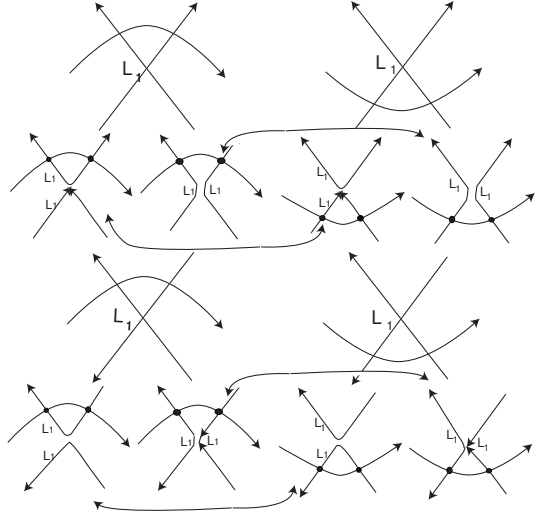} \caption{The
behaviour of orientations under the third Reidemeister move}
\label{Fg}
\end{figure}

In the top picture we see that the  summands from the first pair in
the LHS and RHS contain two crossings with opposite orientations
each, so the total number of crossings contributing to each
orientation is the same for these two summands.

The second in the second pair are identical.

In the bottom picture, the first summand in the LHS has two
crossings contributing the same to the orientation of $L_{1}$; so,
their impact cancels, as well as that for the fist summand of the
RHS. For the second summand in the LHS, we have two crossings with
opposite orientations, and the same in the RHS. So, the orientations
of the corresponding summands in the LHS and in the RHS are the
same.
\end{proof}

Now, we are ready to prove theorem \ref{minexm}. Indeed, if a link
$L=L_{1}\cup L_{2}$ satisfies the conditions of Theorem \ref{minexm}
then $\{L\}_{2}=L$, and if ${\bar L}$ denotes the two-component link
with the orientation of $L_{1}$ reversed then $\{{\bar
L}\}_{2}={\bar {\{L\}_{2}}}={\bar L}$, and since $L$ and ${\bar L}$
are not isomorphic as four-valent framed oriented graphs, they are
not equivalent as framed two component links with one component
oriented.

\begin{figure}
\centering\includegraphics[width=200pt]{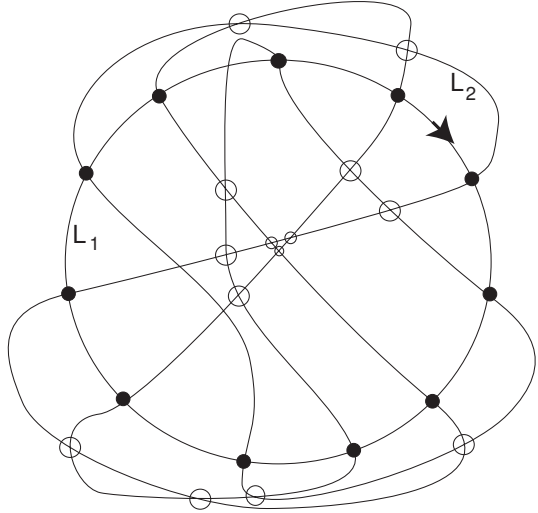} \caption{An
example of a non-orientable free link} \label{exmlink}
\end{figure}

As an example of a link $L$ satisfying the conditions of Theorem
\ref{minexm}, we may consider the link shown in Fig. \ref{exmlink}.

To see that it is indeed non-invertible, let us enumerate the points
on $L_{1}$ along the orientation of $L_{1}$: $A_{1},\dots, A_{11}$
and count whether the distances with respect to $A_{1}$ between
adjacent crossings of $A_{2}$. Thus we get  a sequence of numbers
defined modulo $11$ and up to sign: denote the distance between
$A_{i}$ and $A_{i+1}$ along the component $L_{2}$ by $\beta_{i}$,
$i$ is taken modulo $11$. The numbers $\beta_{i}$ are defined up to
sign because the component $L_{2}$ is not oriented. The (cyclic)
sequence is $3,3,3,4,6,7,6,2,6,9,6$ for one orientation of $L_{2}$.
This sequence has only one fragment of three consequtive equal
numbers: $3,3,3$. If we take the other orientation of $L_{2}$, we
shall get $8,8,8,5,2,5,9,5,4,5,7$ (with three consequtive $8$'s).
None of these two sequences coincides with the cyclic sequences
obtained by inverting the orientation of $L_{1}$: they will have
fragments $6,3,3,3,4$ and $7,8,8,8,5$.

Finally, if we change the roles of $L_{1}$ and $L_{2}$ we shall get
four other cyclic sequences, e.g., $8,6,8,5,4,3,7,7,8,4,4$ (and
similar) which have no three consequtive equal numbers.

So, $L$ (with an orientation of $L_{1}$ fixed) is an example of a
two-component free link with unordered components such that ${\bar
L}$ is not equivalent to $L$.

\section{A Non-Invertible Free Knot}

\begin{figure}
\centering\includegraphics[width=200pt]{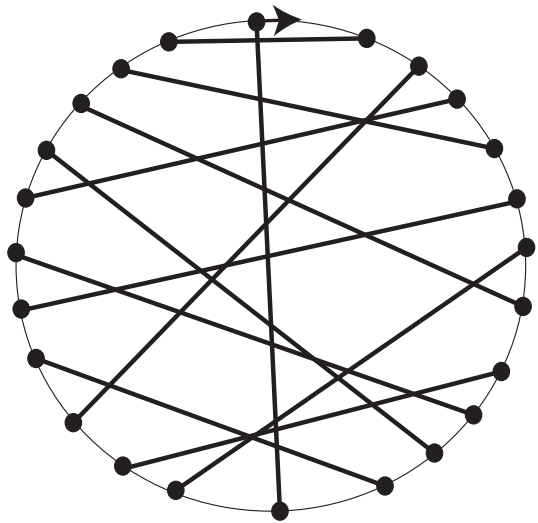} \caption{A
non-invertible free knot}\label{exmknot}
\end{figure}

Consider the Gauss diagram $K$ shown in Fig. \ref{exmknot}.

\begin{st}
The free knot $K$ represented by the diagram shown  in Fig.
\ref{exmknot}, is not equivalent to its inverse.
\end{st}

\begin{proof}
Consider the knot ${\bar K}$ obtained from $K$ by reversion the
orientation.

By construction, we have $\Delta(\bar K)=\overline{ (\Delta(K))}$.

Thus, if we show that $\Delta(K)$ is not invertible as a
$2$-component free link then we see that $K$ is not invertible
either.

Let us extend the map $\{\cdot\}_{2}$ to all two-component free
links. This map is already defined for those links where two
components have an odd intersection. We extend it just by $0$ to the
remaining two-component links.

Now, $\Delta(K)$ is a $\Z_{2}$-linear combination of $2$-component
free links. Moreover, the chord diagram $C(K)$ has exactly one chord
which is linked with all the other chords. The result of smoothing
along this chords leads to the link $L$ from the previous example.

Note, that for this particular $K$ we have
$\Delta(K)=L+\sum_{i}L_{i}$ where $L$ is the link $L$ from the
previous example and all links $L_{i}$ have at least one crossing
belonging to one component.

Now, if we take $\{\Delta(K)\}_{2}$, we get exactly one summand
($L$) which is represented by a diagram with $11$ crossings and can
not be represented by a diagram with a fewer number of crossings,
and all the diagrams $L_{i}$ have strictly less than $11$ crossings.
Analogously, $\{\Delta({\bar K}\}_{2}={\bar L}+\sum {\bar L}_{i}$.

Now, since $L$ is not equivalent to ${\bar L}$ as elements from
$\tLo$ and neither of $L$ or ${\bar L}$ is equivalent to none of
$L_{i}$ or ${\bar L_{i}}$, we see that $\{\Delta({\bar K})\}_{2}\neq
\{\Delta(K)\}_{2}$, so $K$ is not equivalent to ${\bar K}$

\end{proof}

\section{Acknowledgements}

The author is grateful to D.P.Ilyutko, D.Yu.Krylov, and A.Gibson for
fruitful discussions.

\end{document}